\documentclass{article}
\usepackage{amsmath,amsthm,amssymb}
\usepackage[dvipdfmx]{graphicx}
\usepackage{pifont}
\usepackage{wrapfig}
\usepackage{subfigure}

\newtheoremstyle{theorem}
  {10pt}		  
  {10pt}  
  {\sl}  
  {\parindent}     
  {\bf}  
  {. }    
  { }    
  {}     
\theoremstyle{theorem}

\newtheoremstyle{defi}
  {10pt}		  
  {10pt}  
  {\rm}  
  {\parindent}     
  {\bf}  
  {. }    
  { }    
  {}     
\theoremstyle{defi}

\begin{document}

\title{\Large{REVISITING UNIT FRACTIONS THAT SUM TO 1}}

\author{Yutaka Nishiyama\\
Department of Business Information,\\
Faculty of Information Management,\\ 
Osaka University of Economics,\\
2, Osumi Higashiyodogawa Osaka, 533-8533, Japan\\
nishiyama@osaka-ue.ac.jp\\[2pt]}

\maketitle

\noindent
\textbf{Abstract:} This paper is a continuation of a previous paper. Here, as there, we examine the problem of finding the maximum number of terms in a partial sequence of distinct unit fractions larger than 1/100 that sums to 1. In the previous paper, we found that the maximum number of terms is 42 and introduced a method for showing that. In this paper, we demonstrate that there are 27 possible solutions with 42 terms, and discuss how primes show that no 43-term solution exists.



\section*{1. The 27 solutions for n = 42}

Here's a problem I've been thinking about off and on for over twenty years now: Limiting ourselves to distinct positive two-digit denominators, what is the largest number of unit fraction terms whose sum is 1? Using equations, the problem looks like this:

\bigskip
Find the maximal value of $n$ in
$$ 1=\sum^{n}_{i=1} \frac{1}{a_i},  $$

where $a_i$ is a natural number $0<a_1<a_2< \cdots < a_n \le 99.$

\bigskip
For example, in the equation

$$ 1= \frac{1}{2}+\frac{1}{3}+\frac{1}{6}, $$

\noindent
we have $n = 3$.

This problem was first posed in my article \textit{In search of an elegant solution}[1], which was published in the magazine \textit{Mathematics Seminar}. The answers I received from readers were that $n = 42$, up to differences in terms. However, no one provided a proof that no solution with $n = 43$ exists. Since then, I have spent a long time wanting to put this problem to bed by proving this.

I recently took a big step toward solving this problem thanks to a letter I received from Mr. Teruo Nishiyama from Tokyo, who gave 25 possible solutions using 42 terms. His method for showing this begins with separating two-digit numbers that are usable in a 42-term solution from those that are not. This must be determined empirically, classifying as ``usable''  those numbers that appear in solutions created by trial-and-error. He found 50 in all:

$$13, 15, 17, 18, 19, 20, 21, 22, 24, 26,$$
$$27, 28, 30, 32, 33, 34, 35, 36, 38, 39,$$
$$40, 42, 44, 45, 48, 50, 52, 54, 55, 56,$$
$$57, 60, 63, 65, 66, 70, 72, 75, 76, 77,$$
$$78, 80, 84, 85, 88, 90, 91, 95, 96, 99.$$

Next one looks for sums of unit fractions that equal 1, using 43 of the above numbers. When all combinations of unit fractions using 43 of the above numbers are examined, there is no case where their sum is 1, but there are 25 cases where combinations of 42 of these numbers do sum to 1.

The idea of focusing on only those numbers that can be used, rather than all of the numbers, belongs to Mr. Nishiyama, and I did verify that all 25 solutions he proposed give correct sums.

However, another reader, Mr. Shigeru Mori, sent a response to the problem in which he claimed to have found 27 solutions for the $n = 42$ case (though without details regarding the solutions and his method for finding them).

It appears that Mr. Nishiyama had omitted 14 from his list of usable numbers because it had not appeared in previous examples of successful sums. However, as an example, in the case of a 3-term decomposition like
$$ \frac{1}{9}=\frac{1}{14}+\frac{1}{35}+\frac{1}{90} $$

\noindent
a 14 appears on the right side during unit fraction decomposition, and so 14 must be added to the list. This means that there are actually 51 numbers from which to choose, from which we eliminate sets of 9 to search for a 42-term solution.

When I tried this I did indeed find one solution using a 14. I now had 26 solutions of my own, but Mr. Mori claimed 27. Curious as to what the missing solution might be, I began looking for another number to add to the 51 that I knew I could use. One way of finding out would be to take the 90 numbers from 10 to 99 and try exhaustively eliminating 48 numbers from them, but my poor old computer is not up to the task, and besides, as a lover of mathematics I prefer to avoid such brute force methods when possible.

So I was back to the problem of deciding which numbers are usable and which are not. Despite repeated attempts, I found myself unable to come up with a satisfactory method for determining this, so I tried contacting Mr. Koji Oseki of Niigata prefecture, who had contacted me before in regards to this problem. He immediately replied, giving me a solution I had not previously found. It was a solution that contained a 1/12 term, meaning it was 12 I needed to add to my list of usable numbers. (This was particularly interesting because 1/12 is larger than the previously largest value of 1/13.)

So I now had a list of 52 usable numbers, from which 10 should be eliminated to look for 42-term solutions. Doing so I did indeed find a new solution, bringing my tally up to 27, in agreement with Mr. Mori's results.

\section*{2. Usable and unusable numbers}

So what is it about the excluded 47 numbers that make them unusable?

There are 99 numbers from 1 to 99. We can immediately exclude the numbers 1 through 11 as unusable, because the only solution for 1/1 is $n = 1$; additionally, the values of 1/2 through 1/11 are too large to contribute to increasing the number of terms in the unit fraction decomposition. We can therefore limit consideration to integers 12 or greater.

There are 99 numbers in all, among which we have empirically found 52 usable numbers and excluded 11, meaning that there are $99-11-52=36$ more unusable numbers, as follows:
$$16, 23, 25, 29, 31, 37, 41, 43, 46, 47,$$
$$49, 51, 53, 58, 59, 61, 62, 64, 67, 68,$$
$$69, 71, 73, 74, 79, 81, 82, 83, 86, 87,$$
$$89, 92, 93, 94, 97, 98.$$

So what characterizes these numbers? First, we can see that all primes 23 and larger are included in the list. There are 17 of those:
$$23, 29, 31, 37, 41, 43, 47, 53, 59, 61,$$
$$67, 71, 73, 79, 83, 89, 97.$$

The other primes are 13, 17, and 19, each of which appear in unit fractions 1/13, 1/17, 1/19 during the decomposition process, so only primes up to 19 are usable.

Furthermore, we see that the second, third, and fourth multiples of primes $(2p, 3p, 4p)$, although not primes themselves, are also unusable. There are 13 of those:

7 values for $2p$: \;\;\;\;	46, 58, 62, 74, 82, 86, 94

4 values for $3p$: \;\;\;\;	51, 69, 87, 93

2 values for $4p$: \;\;\;\;	68, 92

Many perfect squares are also unusable numbers. There are 5 of these: 
$$16, 25, 49, 64, 81$$

However note that 36 is a perfect square, but is usable. Furthermore, the second multiple of a perfect square is unusable, adding 98 to that list (because $98=2 \times 7^{2} $).

In total, we have found $17+13+5+1=36$ unusable numbers. I will omit here a proof of why primes and multiples of primes are not usable, but I encourage readers to consider this question themselves.

\section*{3. Two- and three-term decompositions }

In this section I would like to present a simple method for decomposition into unit fractions.

Starting with the case of
$$ 1=\frac{1}{2}+\frac{1}{3}+\frac{1}{6}, $$

\noindent
we can decompose the 1/3 and the 1/6 as
$$ \frac{1}{3}=\frac{1}{4}+\frac{1}{12} $$
$$ \frac{1}{6}=\frac{1}{7}+\frac{1}{42}, $$

\noindent
giving us
$$ 1=\frac{1}{2}+\biggl( \frac{1}{4}+\frac{1}{12} \biggr)+\biggl( \frac{1}{7}+\frac{1}{42} \biggr). $$

\noindent
Two- and three-term decompositions like this are basic techniques for looking for longer sums, so let's look into those more closely.

\bigskip
\noindent
\textit{Two-term decompositions}

Say you can decompose the denominator $n$ of a unit fraction into a product as $n=a \times b$ (allowing $a, b$ to be 1). Then, you can decompose the unit fraction into the sum of two unit fractions, as follows:
$$ \frac{1}{n}=\frac{1}{a \times b}=\frac{a+b}{(a \times b)(a+b)}=\frac{1}{b(a+b)}+\frac{1}{a(a+b)} $$

For example, since 3=1~3 we can write this as
$$ \frac{1}{3}=\frac{1}{1 \times 3}=\frac{1+3}{(1 \times 3)(1+3)}=\frac{1}{3(1+3)}+\frac{1}{(1+3)}=\frac{1}{12}+\frac{1}{4}. $$

\bigskip
\noindent
\textit{Three-term decompositions}

Suppose we have a unit fraction $1/n$ that can be decomposed into three unit fractions as follows:
$$ \frac{1}{n}=\frac{1}{a}+\frac{1}{b}+\frac{1}{c}. $$

\noindent
We can use 1/12 as an example:
$$ \frac{1}{12}=\frac{1}{26}+\frac{1}{39}+\frac{1}{52}. $$

However, it can be tricky to find these. The key to doing so is to find three numbers $a, b, c$ whose least common multiple (LCM) is $n$. Using the 1/12 above as an example, the three numbers 3, 4, and 6 have 12 as a LCM. Next we multiply the numerator and denominator of 1/12 by the sum of these numbers $3+4+6=13$. Since 12 is the LCM of 3, 4, 6 the fractions 3/12, 4/12, 6/12 will reduce to lowest terms with a 1 in the numerator.
$$ \frac{1}{12}=\frac{3+4+6}{12 \times (3+4+6)}=\frac{3}{12 \times 13}+\frac{4}{12 \times 13}+\frac{6}{12 \times 13} $$
$$ =\frac{1}{4 \times 13}+\frac{1}{3 \times 13}+ \frac{1}{2 \times 13}=\frac{1}{52}+\frac{1}{39}+\frac{1}{26} $$

We can generalize this as follows:
$$\frac{1}{n}=\frac{a+b+c}{n(a+b+c)}=\frac{a}{n(a+b+c)}+\frac{b}{n(a+b+c)}+\frac{c}{n(a+b+c)}.$$

Here $n$ is the LCM of $a, b, c$, so $n$ is evenly divided by each of $a, b, c$. In other words, the numerator of each of $a/n, b/n, c/n$ is 1. These with denominators multiplied by $(a+b+c)$ are still unit fractions.

There are other methods for finding two- and three-term decompositions. One can also devise ways of finding four- and five-term decompositions, but those will involve the methods for finding two- and three?term decompositions. It is enjoyable to search for such methods using pencil and paper, but to be sure you cover everything itfs probably better to create an exhaustive list in some way, such as by using spreadsheet software.

\section*{4. The relation between the 27 solutions}

So what is the largest number of terms that we might consider possible? As calculated in the previous paper, that number is 62 [2][3], according to the following inequality:
$$ \log n + \frac{1}{n} < \sum^{n}_{i=1} \frac{1}{n} < \log n +1 $$

The larger the denominators in the unit fractions the more terms one can add, so we can consider the case of a sum of unit fractions with terms starting with a denominator of 99 and working backwards, in other words $\displaystyle{\sum^{99}_{i=n} \frac{1}{i}}$. Then we look for the largest value of $n$ such that the sum does not exceed 1. Doing so we find that
$$\frac{1}{38}+\frac{1}{39}+\cdots+\frac{1}{98}+\frac{1}{99}<1,$$

\noindent
but adding $1/37$ makes the inequality false, and from $99-38+1=62$ we see that the upper limit on the number of terms is 62. Not that this is an actual solution, but numerically speaking we can say that more than 62 terms is impossible.

Figure 1 is a listing of all 27 solutions for $n=42$. Looking at this table, we can see that of the 52 usable numbers, 21 are present in all 27 solutions:

$$17, 26, 32, 33, 34, 40, 44, 48, 50, 55,$$
$$56, 66, 75, 76, 77, 80, 84, 85, 88, 91,$$
$$96.$$

Furthermore, 12 and 14 appear in only one solution each; 19 and 57, too, are rarely used, appearing in only three solutions each. 

It is also fun to look for connections between the 27 solutions. For example the first solution is

$$\underline{12}, 17, 21, 22, 24, 26, 27, \underline{30}, 32, 33,$$
$$34, 35, \underline{36}, 38, \underline{39}, 40, 42, 44, 48, 50,$$
$$52, 54, 55, 56, 60, 63, 66, 70, 72, 75,$$
$$76, 77, 78, 80, 84, 85, 88, 90, 91, 95,$$
$$96, 99,$$

\noindent
and the second solution is

$$\underline{13}, 17, \underline{18}, 21, 22, 24, 26, 27, 32, 33,$$
$$34, 35, 38, 40, 42, 44, \underline{45}, 48, 50, 52,$$
$$54, 55, 56, 60, 63, \underline{65}, 66, 70, 72, 75,$$
$$76, 77, 78, 80, 84, 85, 88, 90, 91, 95,$$
$$96, 99.$$

\noindent
The differing numbers (underlined) are

$$12, 30, 36, 39 \;\; \rm{and} \;\; 13, 18, 45, 65,$$

\noindent
and these values have the following relation:
$$ \frac{1}{12}+\frac{1}{30}+\frac{1}{36}+\frac{1}{39}=\frac{1}{13}+\frac{1}{18}+\frac{1}{45}+\frac{1}{65}=\frac{199}{1170}.$$

Each of the 27 solutions are similarly connected through replacement of from two to five of the terms they contain.

\begin{figure}[htbp]
\begin{center}
\includegraphics[width=100mm]{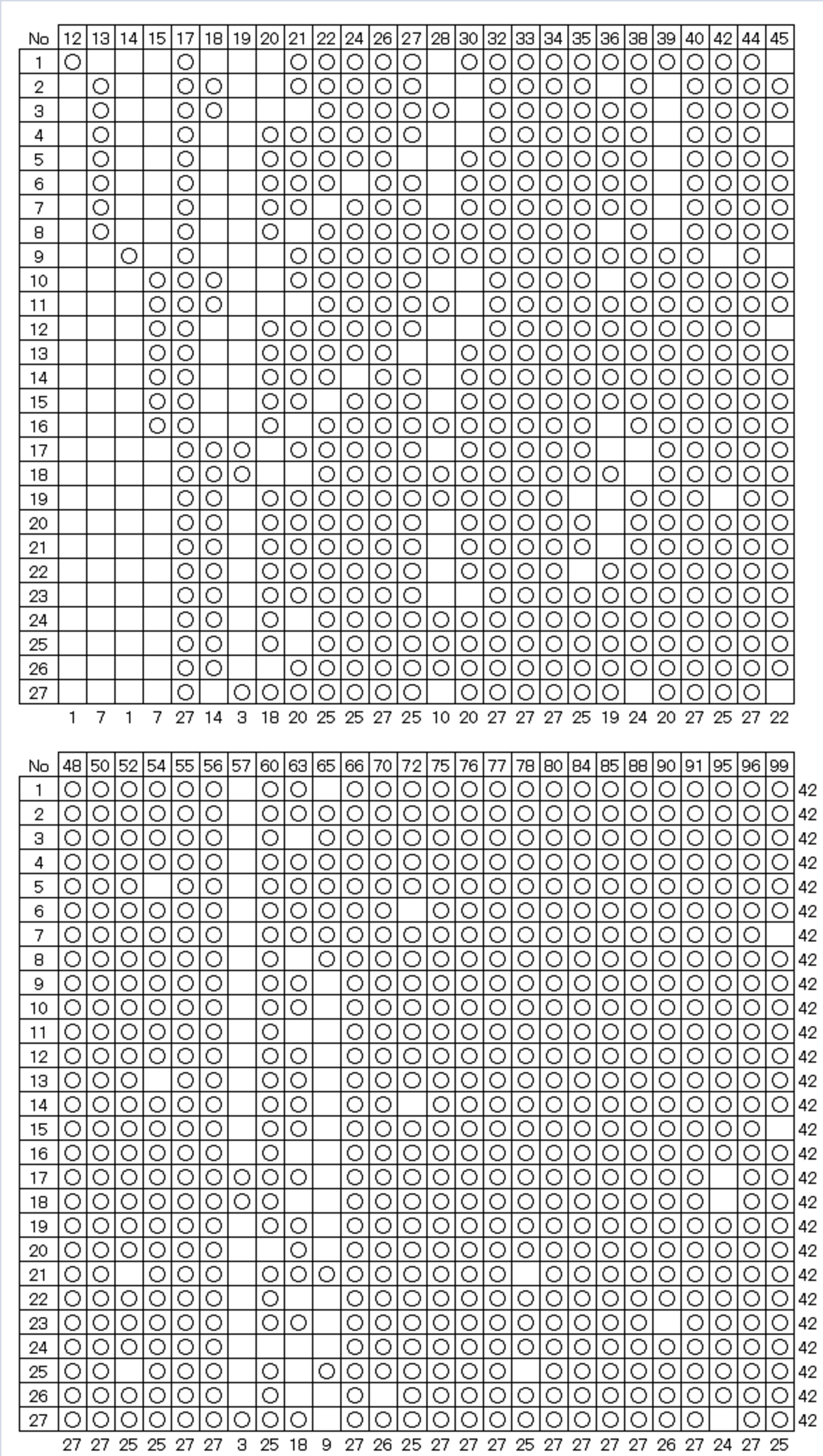}
\end{center}
\caption{The 27 solutions for n = 42.}
\end{figure}

\end{document}